# Weak and strong convergence of a relaxed inertial proximal splitting algorithm for solving hierarchical equilibrium problems*

Zakaria Mazgouri$^{[ORCID:0000-0003-0131-4741]}$,
Hassan Riahi$^{[ORCID:0000-0001-6248-1266]}$, and
Michel Théra$^{[ORCID:0000-0001-9022-6406]}$

**Abstract** In this chapter, we introduce the relaxed inertial proximal splitting algorithm (RIPSA) for hierarchical equilibrium problems. Using Opial-Passty's lemma, we first establish weak ergodic and weak convergence of the sequence generated by the algorithm to a solution of the problem, in the absence of the Browder-Halpern contraction factor. We then derive a strong convergence result under an additional strong monotonicity assumption. Subsequently, we relax this requirement by removing strong monotonicity and instead incorporating a Browder-Halpern contraction factor into (RIPSA), which guarantees strong convergence to a solution determined by the contraction factor. Finally, we discuss two related settings: convex minimization problems and monotone variational inequalities formulated as fixed-point problems for nonexpansive operators.

Zakaria Mazgouri
Laboratory LSATE, National School of Applied Sciences, Sidi Mohamed Ben Abdellah University, 30000, Fez, Morocco
e-mail: zakaria.mazgouri@usmba.ac.ma

Hassan Riahi†Corresponding author)
Laboratory of Mathematics, Modeling and Automatic Systems, Cadi Ayyad university, Faculty of Sciences Semlalia, 40000, Marrakesh, Morocco
e-mail: h-riahi@uca.ac.ma

Michel Théra
XLIM UMR-CNRS 7252, Université de Limoges, Limoges, France
e-mail: michel.thera@unilim.fr

# 1 Introduction

Throughout this paper, $\mathcal{H}$ denotes a real Hilbert space with inner product $\langle \cdot, \cdot \rangle$ and norm $\|\cdot\|$. Let $C \subset \mathcal{H}$ be a nonempty closed convex set, and let $F, G : C \times C \to \mathbb{R}$ be two bifunctions.

We study a two-level hierarchical equilibrium problem in which the upper-level equilibrium is constrained to the solution set of a lower-level problem. The lower equilibrium problem is defined as

$$\text{find } x^* \in C \text{ such that } F(x^*, y) \geq 0, \quad \text{for all } y \in C, \qquad \textbf{(EP)}$$

with solution set denoted by $\mathbf{S}_F$. The hierarchical problem is then formulated as

$$\text{find } \overline{x} \in \mathbf{S}_F \text{ such that } G(\overline{x}, y) \geq 0, \quad \text{for all } y \in \mathbf{S}_F. \qquad \textbf{(HEP)}$$

The theory of equilibrium problems (EPs) has developed along a rich trajectory that naturally connects existence results, stability analysis, dynamical approaches, and iterative algorithms, while at the same time motivating a wide range of applications in optimization, economics, engineering, and network sciences. We briefly recall this trajectory and highlight some key contributions.

**Existence theory.** The origins of EPs trace back to nonlinear analysis. The celebrated Knaster–Kuratowski–Mazurkiewicz (KKM) lemma and Ky Fan's minimax principle [1] provided fundamental existence criteria. Further refinements were given by Brézis, Nirenberg, and Stampacchia [2], and by Mosco [3], who developed variational techniques that remain central today. The unifying contribution of Blum and Oettli [4] marked a turning point, showing that EPs encompass optimization, variational inequalities, fixed-point problems, and Nash equilibria within a single abstract framework.

**Stability and sensitivity analysis.** Beyond existence, stability questions naturally arise: how do solutions behave under perturbations of the underlying data? Early contributions by Ait Mansour and Riahi [5] established sensitivity properties of abstract EPs, ensuring robustness of solutions. Their results can also be interpreted as a pseudo-algorithmic counterpart to iterative procedures. Stability of (**EP**) was later extended to a dynamical setting in [6], while interesting quantitative stability results for the hierarchical problem (**HEP**) were obtained in [7]. More recently, sensitivity analysis has been extended to equilibrium problems with trifunctions [8]. Additional work on monotonicity and generalized monotone bifunctions [9, 10] further strengthened the theoretical foundations and paved the way for reliable computational schemes.

**Dynamical approaches.** An alternative path toward the resolution of EPs lies in continuous-time dynamics. The asymptotic behavior of first-order dynamical systems associated with EPs was studied in [11], establishing convergence of trajectories to equilibrium solutions. This approach was later extended to the hierarchical framework (**HEP**) through the analysis of mixed dynamical equilibrium systems [12].

Such dynamical perspectives not only clarify asymptotic properties but also inspire inertial and time-discretized algorithms.

**Iterative schemes and convergence.** On the algorithmic side, the proximal point method introduced by Martinet [13] and Rockafellar [14] has played a foundational role. Moudafi [15] adapted this approach to equilibrium problems, and later introduced inertial and bilevel variants [16, 17]. Inspired by Antipin and Flåm [18], inertial-type schemes have since become powerful tools to accelerate convergence. A variety of related strategies-including prox-penalization, splitting algorithms, viscosity methods, and forward-backward iterations — have been studied in both Hilbert and Banach spaces [19–21]. Convergence analyses often rely on classical fixed-point principles such as Opial's lemma [22], and the works of Browder [23] and Halpern [24]. Recent advances refine the trade-off between weak and strong convergence, stability, and computational efficiency [25–27].

**Applications.** Equilibrium models appear in a broad range of applied sciences. In economics and game theory, they underpin Nash equilibria and market equilibria. In optimization and control, bilevel and hierarchical EPs capture leader-follower dynamics. In mechanics and engineering, they describe models of contact, friction, and network flows. A particularly active domain is electricity markets, where EPs model pricing mechanisms, demand elasticity, and transmission constraints. For example, Aussel and Riccardi [28] applied equilibrium and complementarity formulations to pay-as-clear electricity markets, incorporating demand response and stability considerations. Transportation and traffic networks provide another classical application, where EPs describe user equilibria in congested systems, route choice problems, and network design, linking variational inequalities with practical planning models [29]. Comprehensive treatments of both theory and applications can be found in monographs such as Yuan [30], Göpfert et al. [31], and Bigi et al. [32]. Thus, EPs provide a unifying bridge between abstract variational analysis and concrete applied models.

To summarize, the development of equilibrium problems follows a natural trajectory: existence results guarantee solvability, stability theory ensures robustness, dynamical systems connect continuous-time and discrete-time perspectives, iterative methods enable computation, and applications demonstrate practical relevance. This framework motivates our forthcoming contributions on inertial proximal methods for hierarchical equilibrium problems.

As mentioned earlier, the convergence of iterative schemes for solving (**EP**) has been extensively investigated (see, e.g., [4, 9, 30, 33]). A classical tool in this direction is the proximal method. Moudafi [15] established weak convergence to a solution of (**EP**) via the proximal point algorithm

$$r_k F(x_{k+1}, y) + \langle x_{k+1} - x_k, y - x_{k+1} \rangle \geq 0, \quad \forall y \in C,$$

which naturally induces the resolvent operator

$$J_r^F(x) := \{z \in C : rF(z, y) + \langle y - z, z - x \rangle \geq 0 \ \forall y \in C\}, \quad x \in \mathcal{H}, r > 0,$$

leading to the iteration $x_{k+1} = J_r^F(x_k)$. Later, inspired by Antipin [18], Moudafi [16] introduced the inertial Mann-type iteration

$$x_{k+1} = J_{r_k}^F\big(x_k + \alpha_k(x_k - x_{k-1})\big).$$

Building on these developments, several algorithms were introduced for hierarchical problems. Moudafi [17] proposed the <u>regularized proximal point algorithm</u>

$$x_{k+1} = J_{r_k}^{F+\mu_k G}(x_k),$$

which converges weakly under suitable summability conditions. Chbani and Riahi [25] obtained weak convergence for the variant

$$x_{k+1} = J_{r_k}^{\lambda_k F + G}(x_k),$$

under weaker assumptions and a geometric condition on $F$, and later proved strong convergence when $G$ is strongly monotone. Splitting-type methods were also studied, such as

$$x_{k+1} = J_{r_k}^{\lambda_k F} \circ J_{r_k}^G(x_k),$$

and inertial regularized variants [26, 27], which achieve weak or strong convergence under conditions on $\{r_k\}, \{\lambda_k\}, \{\tau_k\}$:

$$\liminf_{k \to +\infty} r_k > 0, \quad \sum_{k=0}^{+\infty} r_k \mu_n < +\infty, \quad \text{and} \quad \|x_{k+1} - x_k\| = o(\mu_k).$$

Chbani and Riahi [25] proposed a variant that interchanges the penalty parameter between $F$ and $G$, obtaining weak convergence under milder assumptions. Their *regularized proximal point algorithm*

$$x_{k+1} = J_{r_k}^{\lambda_k F + G}(x_k)$$

converges weakly to a solution of (**HEP**), provided that $\liminf_{k \to +\infty} r_k > 0$, $\liminf_{k \to +\infty} \lambda_k = +\infty$, and the following geometric condition on $F$ holds: for every $x \in \mathbf{S}_F$ and $p \in \mathcal{N}_{\mathbf{S}_F}(x)$,

$$\sum_{k=1}^{+\infty} r_k \lambda_k \left[ \mathcal{F}_F\left(x, \tfrac{2}{\lambda_k}p\right) - \sigma_{\mathbf{S}_F}\left(\tfrac{2}{\lambda_k}p\right) \right] < +\infty, \tag{1}$$

where

$$\mathcal{F}_F(x, q) := \sup_{v \in C} \{\langle q, v \rangle + F(v, x)\},$$

and $\mathcal{N}_{\mathbf{S}_F}(x)$ denotes the normal cone to $\mathbf{S}_F$ at $x$, and $\sigma_{\mathbf{S}_F}$ is the support function of $\mathbf{S}_F$. Without assuming condition (1), they further proved that if $G$ is strongly monotone, then $\{x_k\}$ converges strongly to the unique solution of (**HEP**) whenever $r_k \to 0$, $\lambda_k \to +\infty$, $\sum r_k = +\infty$, and $\liminf_k r_k \lambda_k > 0$.

A more efficient scheme considered in that work is the splitting proximal algorithm

$$x_{k+1} = J_{r_k}^{\lambda_k F} \circ J_{r_k}^{G}(x_k),$$

which converges weakly under condition (1) together with

$$\sum r_k = +\infty, \quad \sum r_k^2 < +\infty, \quad \text{and} \quad \liminf_k r_k \lambda_k > 0.$$

Strong convergence is also ensured if $r_k \to 0$, $\sum r_k = +\infty$, $\liminf_k r_k \lambda_k > 0$, and condition (1) is replaced by the weaker requirement: for every $x \in \mathbf{S}_F$ and $p \in \mathcal{N}_{\mathbf{S}_F}(x)$,

$$\lim_{k \to +\infty} r_k \lambda_k \left[ \mathcal{F}_F\left(x, \tfrac{2}{\lambda_k} p\right) - \sigma_{\mathbf{S}_F}\left(\tfrac{2}{\lambda_k} p\right) \right] = 0.$$

Another important development incorporates inertial terms into the regularized proximal point algorithm, leading to the iteration

$$x_{k+1} = J_{r_k}^{\lambda_k F + G}\left(x_k + \tau_k(x_k - x_{k-1})\right), \quad \tau_k \in [0,1].$$

Under condition (1), $\liminf_k r_k > 0$, and $\lambda_k \to +\infty$, weak convergence was established in [26, Theorem 3.3]. When $G$ is strongly monotone, strong convergence to the unique solution of (**HEP**) follows under the same conditions, provided $\tau_k < 1/3$. Balhag et al. [27] proposed a modification of this inertial scheme:

$$x_{k+1} = J_{r_k}^{\lambda_k F} \circ J_{r_k}^{G}\left(x_k + \tau_k(x_k - x_{k-1})\right),$$

achieving similar weak and strong convergence results, under conditions analogous to those in [26], together with the additional requirement that $\{\tau_k\}$ is nondecreasing and bounded above by $(\sqrt{3}-1)/4$.

Motivated by these contributions, this contribution develops a new scheme ensuring strong convergence to a solution of (**HEP**). Unlike [27], where strong monotonicity of $G$ is required, we assume only monotonicity and incorporate a relaxation factor into the algorithm. Our method, termed the Relaxed Inertial Proximal Algorithm (RIPA), is defined by

$$x_{k+1} = \alpha_k x_k + \beta_k J_{r_k}^{\lambda_k F} \circ J_{r_k}^{G}\left(x_k - \gamma_k(x_k - x_{k-1})\right) + \mu_k g(x_k),$$

where $r_k > 0$, $\alpha_k, \beta_k, \mu_k \geq 0$, and $\alpha_k + \beta_k + \mu_k = 1$.

**Main contributions.**

- In Section 2, we review preliminary concepts and basic results on monotone bifunctions, together with some technical lemmas.
- Section 3 introduces the Relaxed Inertial Proximal Algorithm (RIPA) and establishes a key control lemma.
- Section 4 provides weak convergence results for the sequence generated by (RIPA) in the case $\alpha_k = 0, \beta_k = 1$.

- Section 5 establishes strong convergence: first when $G$ is $\rho$-strongly monotone and $\alpha_k = 0, \beta_k = 1$, and then when $\alpha_k, \beta_k, \mu_k$ are positive, $g$ is a $\iota$-contraction with $0 < \iota < 1$, and $G$ is only monotone (i.e., $\rho = 0$).
- Section 6 discusses applications to convex minimization, variational inequalities, and fixed point problems.

## 2 Preliminaries

### 2.1 Basic definitions

**Subdifferential.** Let $f : \mathcal{H} \to \mathcal{H}$ be a convex function. The underline{subdifferential} of $f$ at $x \in \mathcal{H}$ is defined as

$$\partial f(x) := \{ y \in \mathcal{H} : f(z) \geq f(x) + \langle y, z - x \rangle \text{ for all } z \in \mathcal{H} \}.$$

**Support function.** The support function of a nonempty set $C \subseteq \mathcal{H}$ is given by

$$\sigma_C(y) := \sup_{x \in C} \langle y, x \rangle, \quad y \in \mathcal{H}.$$

**Normal cone.** Let $C \subseteq \mathcal{H}$ be a nonempty closed convex set, and let $\bar{x} \in C$. A vector $d \in \mathcal{H}$ is called a normal vector to $C$ at $\bar{x}$ if

$$\langle d, x - \bar{x} \rangle \leq 0 \quad \text{for all } x \in C.$$

The set of all such vectors, called the normal cone to $C$ at $\bar{x}$, is denoted by $\mathcal{N}_C(\bar{x})$ and is a convex cone:

$$\mathcal{N}_C(\bar{x}) := \{ d \in \mathcal{H} : \langle d, x - \bar{x} \rangle \leq 0 \text{ for all } x \in C \}.$$

**Metric projection.** For a nonempty closed convex set $C \subseteq \mathcal{H}$, the metric projection $\Pi_C : \mathcal{H} \to C$ is defined as

$$\Pi_C(x) := \arg\min_{y \in C} \|x - y\|.$$

It is well known that $\Pi_C$ is single-valued and characterized by

$$y = \Pi_C(x) \iff y \in C \text{ and } \langle x - y, z - y \rangle \leq 0 \quad \text{for all } z \in C,$$

which is equivalent to $x - y \in \mathcal{N}_C(y)$.

Moreover, the projection operator $\Pi_C$ is firmly nonexpansive, i.e.,

$$\|\Pi_C(x) - \Pi_C(y)\|^2 \leq \langle \Pi_C(x) - \Pi_C(y), x - y \rangle \quad \text{for all } x, y \in \mathcal{H},$$

and hence also nonexpansive:

$$\|\Pi_C(x) - \Pi_C(y)\| \leq \|x - y\|, \quad \text{for all } x, y \in \mathcal{H}.$$

**Indicator function.** Let $C \subseteq \mathcal{H}$. The underline{indicator function} of $C$, denoted $\iota_C : \mathcal{H} \to \mathbb{R} \cup \{+\infty\}$, is defined by

$$\iota_C(x) := \begin{cases} 0, & \text{if } x \in C, \\ +\infty, & \text{if } x \notin C. \end{cases}$$

**Fenchel conjugate.** Let $f : \mathcal{H} \to \mathbb{R} \cup +\{\infty\}$ be a convex function. The Fenchel conjugate $f^*$ of $f$ is defined by

$$f^*(x) = \sup_{y \in \mathcal{H}} \{\langle y, x \rangle - f(y)\}.$$

**Resolvents of monotone operators.** Let $A : \mathcal{H} \rightrightarrows \mathcal{H}$ be a (possibly multivalued) monotone operator on a Hilbert space $\mathcal{H}$. For $\lambda > 0$, the underline{resolvent} of $A$ with parameter $\lambda$ is defined by

$$J_\lambda^A := (I + \lambda A)^{-1}.$$

If $A$ is maximally monotone, then $J_{\lambda A}$ is single-valued and firmly nonexpansive. In particular, when $A = \mathcal{N}_C$ is the normal cone to a closed convex set $C$, the resolvent coincides with the metric projection:

$$J_\lambda^{\mathcal{N}_C} = \Pi_C, \qquad \text{for all } \lambda > 0.$$

**Proximal operators.** A fundamental example of resolvents arises from convex analysis: if $A = \partial f$ is the subdifferential of a proper, convex, lower semicontinuous function $f : \mathbb{R}^n \to (-\infty, +\infty]$, then the resolvent of $A$ is the underline{proximal operator} of $f$:

$$\text{prox}_{\lambda f}(x) := J_\lambda^{\partial f}(x) = \arg\min_{y \in \mathcal{H}} \left\{ f(y) + \tfrac{1}{2\lambda} \|y - x\|^2 \right\}.$$

Thus, the proximal operator generalizes the metric projection (which corresponds to the indicator function of a convex set).

**Remark.** These definitions are not merely abstract: projections, resolvents, and proximal operators constitute the main building blocks of the iterative algorithms studied later in this chapter. In particular, they provide the analytical tools for constructing and analyzing algorithms designed to solve equilibrium problems (**EP**) and hierarchical equilibrium problems (**HEP**).

## 2.2 Basic assumptions on bifunctions

Let $C \subset \mathcal{H}$ be a nonempty closed convex set, and let $F : C \times C \to \mathbb{R}$ be a bifunction. We consider the following conditions:

(C1) $F(x,x) = 0$ for all $x \in C$;
(C2) for each $x \in C$, the mapping $y \mapsto F(x,y)$ is convex;
(C3) for each $x \in C$, the mapping $y \mapsto F(x,y)$ is lower semicontinuous;
(C4) $F$ is monotone, i.e.,

$$F(x,y) + F(y,x) \leq 0, \quad \text{for all } y \in C;$$

(C5) $F$ is upper hemicontinuous, i.e.,

$$\lim_{t \downarrow 0} F(tz + (1-t)x, y) \leq F(x,y), \quad \text{for all } x, y, z \in C;$$

(C6) $F$ is $\rho$-monotone ($\rho \in \mathbb{R}$), i.e.,

$$F(x,y) + F(y,x) \leq -\frac{\rho}{2}\|x-y\|^2, \quad \text{for all } x, y \in C.$$

(C7) $F$ is maximal, i.e., whenever $(x,u) \in C \times \mathcal{H}$ and $F(x,y) \leq \langle u, x-y \rangle$ for all $y \in C$, it follows that $F(x,y) + \langle u, x-y \rangle \geq 0$ for all $y \in C$.

*Remark 1* If $\rho > 0$, then $F$ is called $\rho$-strongly monotone; strong monotonicity implies monotonicity, corresponding to $\rho = 0$. If $\rho < 0$, we say that $F$ is $\rho$-weakly monotone.

### 2.3 Minty's lemma and maximality

The following equivalence is classical.

**Lemma 1 (Minty's lemma)** *[4] Let $C$ be a nonempty closed convex subset of $\mathcal{H}$ and let $F : C \times C \to \mathbb{R}$. Consider*

$$(EP): \text{ find } \bar{x} \in C \text{ such that } F(\bar{x}, y) \geq 0, \text{ for all } y \in C,$$
$$(DKFI): \text{ find } \bar{x} \in C \text{ such that } F(y, \bar{x}) \leq 0, \text{ for all } y \in C.$$

1. *If $F$ satisfies (C4), then every solution of (EP) is also a solution of (DKFI).*
2. *If $F$ satisfies (C1), (C2), and (C5), then every solution of (DKFI) is also a solution of (EP).*

**Lemma 2 ([34])** *Let $C$ be a nonempty closed convex subset of $\mathcal{H}$ and let $F : C \times C \to \mathbb{R}$ satisfy (C1)–(C4). Then the following are equivalent:*

1. *$F$ is maximal;*
2. *for each $x \in \mathcal{H}$ and $r > 0$, there exists a unique $z = J_r^F(x) \in C$ such that*

$$rF(z,y) + \langle y-z, z-x \rangle \geq 0, \quad \text{for all } y \in C.$$

It is known (see [9, Lemma 2.1]) that maximality of $F$ holds under (C1), (C2), and (C5). In this case the resolvent mapping $J_r^F : \mathcal{H} \to C$ is well defined, and

$$\bar{x} = J_r^F(\bar{x}) \iff F(\bar{x}, y) \geq 0, \text{ for all } y \in C.$$

### 2.4 Yosida approximation

The Yosida approximation associated to $F$ is defined by

$$A_r^F(x) = \frac{1}{r}\left(x - J_r^F(x)\right) = \left(rI + (A^F)^{-1}\right)^{-1}(x),$$

where the operator $A^F$ is given by

$$A^F(x) = \begin{cases} \{z \in \mathcal{H} : F(x, y) + \langle z, x - y \rangle \geq 0, \text{ for all } y \in C\}, & x \in C, \\ \emptyset, & x \notin C. \end{cases}$$

The operator $A^F$ has closed convex values, and $A^F$ is monotone if and only if $F$ is monotone, provided $\text{dom}(A^F) = C$.

**Lemma 3** ([11, Lemma 2.5], [19]) *Assume $F$ satisfies (C1)–(C4). Define $A_o^F(x) := \Pi_{A^F(x)}(0)$, the projection of $0$ onto $A^F(x)$. Then:*
*(i) $\|A_r^F(x)\| \leq \|A_o^F(x)\|$, for all $x \in C$ and $r > 0$;*
*(ii) $r\|A_r^F(x) - A_r^F(y)\|^2 \leq \langle A_r^F(x) - A_r^F(y), y - x \rangle$, for all $x, y \in C$ and $r > 0$;*
*(iii) $\|J_r^F(x) - J_r^F(y)\| \leq \|x - y\|$, for all $x, y \in \mathcal{H}$ and $r > 0$.*

### 2.5 Auxiliary lemmas

We also recall several technical results frequently used in convergence analysis.

**Lemma 4** *For all $x, y, z \in \mathcal{H}$ and $\alpha, \beta, \gamma \in \mathbb{R}$ with $\alpha + \beta + \gamma = 1$,*

$$\|\alpha x + \beta y + \gamma z\|^2 = \alpha\|x\|^2 + \beta\|y\|^2 + \gamma\|z\|^2 - \alpha\beta\|x-y\|^2 - \beta\gamma\|y-z\|^2 - \alpha\gamma\|x-z\|^2.$$

**Lemma 5** *Let $\{u_k\}$ and $\{v_k\}$ be real sequences. If $\sum u_k v_k < +\infty$ and $\sum u_k = +\infty$, then*

$$\liminf_{k \to \infty} v_k \leq 0.$$

**Lemma 6** ([20]) *Let $\{\alpha_k\} \subset (0, 1)$ satisfy $0 < \liminf \alpha_k \leq \limsup \alpha_k < 1$. Let $\{x_k\}$ and $\{v_k\}$ be bounded sequences in $\mathcal{H}$ such that*

$$x_{k+1} = \alpha_k x_k + (1 - \alpha_k) v_k,$$

*and*

$$\limsup_{k \to \infty} \left(\|v_{k+1} - v_k\| - \|x_{k+1} - x_k\|\right) \leq 0.$$

*Then $v_k - x_k \to 0$ strongly.*

**Lemma 7 ([35, Lemma 5])** *Let $\{a_k\}, \{b_k\}, \{\delta_k\}, \{\theta_k\}$ be nonnegative sequences with $\delta_k \in [0, 1/2]$ and $\sum_{k=k_0}^{\infty} \theta_k < +\infty$. Suppose for all $k \geq k_0 > 0$,*

$$a_{k+1} + b_k \leq (1 - \delta_k)a_k + \delta_k a_{k-1} + \theta_k.$$

*Then:*

1. *$\{a_k\}$ converges and $\sum b_k < +\infty$;*
2. *if $\liminf a_k = 0$, then $a_k \to 0$. A sufficient condition is the existence of a positive sequence $\{\lambda_k\}$ such that $\sum \lambda_k a_k < \infty$ and $\sum \lambda_k = \infty$.*

**Lemma 8 ([35, Lemma 6])** *Let $\{a_k\}, \{b_k\}, \{\delta_k\}, \{\theta_k\}, \{s_k\}, \{t_k\}$ be sequences of nonnegative real numbers satisfying*

$$\delta_k \in [0, 1/2], \quad \limsup_{k \to \infty} s_k \leq 0, \quad \sum \theta_k < +\infty, \quad \sum t_k = +\infty, \qquad (2)$$

*and for each $k \geq k_0$,*

$$a_{k+1} + b_k \leq (1 - t_k - \delta_k)a_k + \delta_k a_{k-1} + t_k s_k + \theta_k.$$

*Then $a_k \to 0$ and $\sum b_k = +\infty$.*

We say that a sequence $\{x_k\} \subset \mathcal{H}$ converges <u>weakly ergodically</u> if there exists a nonsummable sequence of positive real numbers $\{r_k\}$ such that the associated ergodic average

$$\tilde{x}_n := \frac{1}{\sigma_n} \sum_{k=1}^{n} r_k x_k, \quad \text{where } \sigma_n = \sum_{k=1}^{n} r_k,$$

converges weakly in $\mathcal{H}$.

We now combine two classical results - Passty's weak ergodic convergence lemma and Opial's weak convergence lemma - into the following statement.

**Lemma 9 (Opial-Passty Lemma [22, 36])** *Let $\{x_k\}$ be a sequence in a Hilbert space $\mathcal{H}$, and let $B \subset \mathcal{H}$ be nonempty. Suppose that:*

1. *for every $y \in B$, the sequence $\{\|x_k - y\|\}$ converges;*
2. *every weak cluster point of $\{x_k\}$ (respectively of $\{\tilde{x}_n\}$) belongs to $B$.*

*Then $\{x_k\}$ (respectively $\{\tilde{x}_n\}$) converges weakly to some $\overline{x} \in B$.*

## 3 Relaxed inertial proximal algorithm for hierarchical equilibrium problems

For the hierarchical equilibrium problem (**HEP**), we will use, depending on the context some of the following assumptions:

(*H*1) $F$ and $G$ satisfy the conditions $(C1) - (C4)$;
(*H*2) for each $x \in C$, $\partial G_x(x) \neq \emptyset$, where $G_x(u) = G(x, u)$, for all $u \in C$ and $+\infty$ if $u \notin C$;
(*H*3) $\mathbb{R}_+(C - \mathbf{S}_F) := \bigcup_{\varepsilon > 0} \varepsilon(C - \mathbf{S}_F)$ is a closed linear subspace of $\mathcal{H}$;
(*H*4) for each $x \in \mathbf{S}, p \in \mathcal{N}_{\mathbf{S}_F}(x)$,

$$\sum_{k=1}^{+\infty} \lambda_k r_k \left[ \mathcal{F}_F\left(x, \frac{2p}{\lambda_k}\right) - \sigma_{\mathbf{S}_F}\left(\frac{2p}{\lambda_k}\right) \right] < +\infty;$$

(*H*5) for each $x \in \mathbf{S}, p \in \mathcal{N}_{\mathbf{S}_F}(x)$,

$$\limsup_{k \to +\infty} \lambda_k r_k \left[ \mathcal{F}_F\left(x, \frac{2p}{\lambda_k}\right) - \sigma_{\mathbf{S}_F}\left(\frac{2p}{\lambda_k}\right) \right] \leq 0.$$

*Remark 2* 1. The function

$$\mathcal{F}_F(x, p) := \sup_{v \in C} (\langle p, v \rangle + F(v, x))$$

is the *Fitzpatrick transform* of $F$. If $F$ satisfies $(H1)$ and we define

$$\bar{F}(x, v) := \begin{cases} F(x, v), & v \in C, \\ +\infty, & v \notin C, \end{cases}$$

then

$$\mathcal{F}_F(x, p) \leq \sup_{v \in C}(\langle p, v \rangle - F(x, v)) = \sup_{v \in \mathcal{H}}(\langle p, v \rangle - \bar{F}(x, v)) = \bar{F}(x, \cdot)^*(p),$$

where $\bar{F}(x, \cdot)^*$ denotes the Fenchel conjugate of $\bar{F}(x, \cdot)$.
2. Assumption $(H3)$ ensures the qualification condition

$$\partial(G_x + \iota_{\mathbf{S}_F})(x) = \partial G_x(x) + \mathcal{N}_{\mathbf{S}_F}(x), \quad x \in \mathbf{S},$$

see [37, 38].
3. Section 6 discusses assumptions $(H4)$-$(H5)$ for the special cases where the lower-level problem is a minimization problem or a fixed point problem.

We now present the **R**elaxed **I**nertial **P**roximal **A**lgorithm (RIPA) for solving (**HEP**). The method depends on a sequence $\{r_k\}$ of positive parameters and sequences $\{\gamma_k\}, \{\alpha_k\}, \{\beta_k\} \subset [0, 1]$.

**Algorithm (RIPA)**

(1i) Initialize $x_0, x_1 \in C$.
(2i) Set $y_k = (1 - \gamma_k)x_k + \gamma_k x_{k-1}$.
(3i) Compute $u_k = J_{r_k}^G(y_k)$ and $z_k = J_{r_k}^{\lambda_k F}(u_k)$.

(4i) Update
$$x_{k+1} = \alpha_k x_k + \beta_k z_k + (1 - \alpha_k - \beta_k) g(x_k),$$

where $g : C \to C$ is a contraction.

Equivalently, the iteration can be written compactly as

$$x_{k+1} = \alpha_k x_k + \beta_k J_{r_k}^{\lambda_k F} \circ J_{r_k}^G \big((1 - \gamma_k)x_k + \gamma_k x_{k-1}\big) + (1 - \alpha_k - \beta_k) g(x_k).$$

*Remark 3* By Lemma 2, the transition $(x_{k-1}, x_k) \mapsto x_{k+1}$ is well-defined for all $r_k > 0$ whenever $\alpha_k, \gamma_k \geq 0$ and $\beta_k > 0$, assuming that $F$ satisfies $(C1)$–$(C4)$.

Let us mention some special cases covered by (RIPA):

- If $F \equiv 0$, then $\mathbf{S}_F = C$ and $J_{r_k}^{\lambda_k F}(u) = \Pi_C(u)$ for every $u \in \mathcal{H}$. For $\beta_k = 1$, $\alpha_k = 0$, and $\gamma_k = 0$, we have $y_k = x_k$ and $u_k = J_{r_k}^G(x_k) \in C$. Hence

$$z_k = J_{r_k}^{\lambda_k F}(u_k) = u_k = J_{r_k}^G(x_k),$$

so that (RIPA) reduces to
$$x_{k+1} = J_{r_k}^G(x_k),$$

which is precisely the proximal algorithm for equilibrium problems proposed in [15].

- If $1 - \alpha_k - \beta_k = 0$ and $\gamma_k = 0$, then (RIPA) reduces to

$$x_{k+1} = \alpha_k x_k + (1 - \alpha_k) J_{r_k}^F(x_k),$$

which is a generalized Mann-proximal algorithm for equilibrium problems.

- If $\rho \in [0, 1)$ and $\gamma_k = 0$, then (RIPA) becomes

$$x_{k+1} = \alpha_k x_k + \beta_k J_{r_k}^F(x_k) + (1 - \alpha_k - \beta_k) g(x_k),$$

which corresponds to the algorithm studied in [21]. In the particular case $g \equiv u$ for some fixed $u \in C$ (i.e., $\rho = 0$), this reduces to the classical Browder-Halpern iteration for fixed points [23, 24].

- For two-level hierarchical equilibrium problems, (RIPA) encompasses several well-known schemes from the literature:

$$\begin{aligned}
x_{k+1} &= J_{r_k}^{F + \mu_k G}(x_k) & [17], \\
x_{k+1} &= J_{r_k}^{\lambda_k F + G}(x_k) & [25], \\
x_{k+1} &= J_{r_k}^{\lambda_k F} \circ J_{r_k}^G(x_k) & [25], \\
x_{k+1} &= J_{r_k}^{\lambda_k F + G}(x_k + \tau_k(x_k - x_{k-1})) & [26], \\
x_{k+1} &= J_{r_k}^{\lambda_k F} \circ J_{r_n}^G(x_k + \tau_k(x_k - x_{k-1})) & [27].
\end{aligned}$$

For the sake of simplifying the convergence proof of (RIPA), we first establish the following control lemma.

**Lemma 10** *Let $F, G : C \times C \to \mathbb{R}$ satisfy (H1). Let $\{x_k\}, \{y_k\}, \{u_k\}, \{z_k\}$ be the sequences generated by algorithm (RIPSA). Assume that the solution set $\mathbf{S}$ of (HEP) is nonempty, and that (H2), (H3) hold. If $G$ is $\rho$-monotone ($\rho \geq 0$), then for each $\{b_k\} \subset (\gamma_k, 1)$ and for every $x \in \mathbf{S}_F$, $v \in A^G(x) + \mathcal{N}_{\mathbf{S}_F}(x)$ and $p \in \mathcal{N}_{\mathbf{S}_F}(x)$ such that $v - p \in A^G(x)$, we have*

$$(1+2\rho r_k)\|z_k - x\|^2 \leq r_k \lambda_k F(z_k, x) + r_k \lambda_k \left[ \mathcal{F}_F\left(x, \tfrac{2}{\lambda_k} p\right) - \sigma_{\mathbf{S}_F}\left(\tfrac{2}{\lambda_k} p\right) \right]$$
$$+ \tfrac{2}{c} r_k^2 \|p\|^2 + (1-\gamma_k)\|x_k - x\|^2 + \gamma_k \|x_{k-1} - x\|^2$$
$$- \left(1 - \tfrac{\gamma_k}{b_k}\right)\|u_k - x_k\|^2 - \gamma_k(1-b_k)\|x_k - x_{k-1}\|^2$$
$$- \left(1 - \tfrac{c}{2}\right)\|u_k - z_k\|^2 + 2\langle v, x - x_k \rangle, \qquad (3)$$

*and*

$$(1+2\rho r_k)\|x_{k+1} - x\|^2 \leq \beta_k r_k \lambda_k F(z_k, x) + \beta_k r_k \lambda_k \left[ \mathcal{F}_F\left(x, \tfrac{2}{\lambda_k} p\right) - \sigma_{\mathbf{S}_F}\left(\tfrac{2}{\lambda_k} p\right) \right]$$
$$+ \tfrac{2}{c} r_k^2 \beta_k \|p\|^2 + \left(\alpha_k(1+2\rho r_k) + \beta_k(1-\gamma_k)\right)\|x_k - x\|^2$$
$$+ \beta_k \gamma_k \|x_{k-1} - x\|^2 - \beta_k \gamma_k (1-b_k)\|x_k - x_{k-1}\|^2$$
$$- \beta_k \left(1 - \tfrac{c}{2}\right)\|u_k - z_k\|^2 - \beta_k \left(1 - \tfrac{\gamma_k}{b_k}\right)\|u_k - x_k\|^2$$
$$+ (1+2\rho r_k)(1-\alpha_k - \beta_k)\|g(x_k) - x\|^2$$
$$- (1+2\rho r_k)\alpha_k \beta_k \|z_k - x_k\|^2 + \frac{2r_k \beta_k}{1+2\rho r_k} \langle v, x - u_k \rangle. \qquad (4)$$

**Proof. 1.** To prove (3), let $x \in \mathbf{S}_F$, $v \in A^G(x) + \mathcal{N}_{\mathbf{S}_F}(x)$ and $p \in \mathcal{N}_{\mathbf{S}_F}(x)$ such that $v - p \in A^G(x)$. Then, for every $k \geq 1$,

$$r_k G(x, y) + r_k \langle v - p, x - y \rangle \geq 0, \quad \text{for all } y \in \mathbf{C}. \qquad (5)$$

The iteration $u_k = J^G_{r_k}(y_k)$ gives

$$r_k G(u_k, y) + \langle u_k - y_k, y - u_k \rangle \geq 0, \quad \text{for all } y \in \mathbf{C}. \qquad (6)$$

Taking $y = u_k$ in (5) and $y = x$ in (6), and summing the inequalities, the $\rho$-monotonicity of $G$ yields

$$-\rho r_k \|x - u_k\|^2 + r_k \langle v - p, x - u_k \rangle + \langle u_k - y_k, x - z_k \rangle + \langle u_k - y_k, z_k - u_k \rangle \geq 0.$$

From (3i) in (RIPSA), we also have

$$r_k \lambda_k F(z_k, x) + \langle z_k - u_k, x - z_k \rangle \geq 0.$$

Summing these two inequalities gives

$$0 \le r_k \lambda_k F(z_k,x) + \langle z_k - y_k, x - z_k \rangle - \rho r_k \|x - u_k\|^2 \\ + r_k \langle -p, x - u_k \rangle + \langle u_k - y_k, z_k - u_k \rangle + r_k \langle v, x - u_k \rangle. \qquad (7)$$

Using the identity $2\langle a,b \rangle = \|a+b\|^2 - \|a\|^2 - \|b\|^2$, we compute

$$2\langle z_k - y_k, x - z_k \rangle = \|y_k - x\|^2 - \|z_k - y_k\|^2 - \|z_k - x\|^2,$$

and

$$2\langle u_k - y_k, z_k - u_k \rangle = \|z_k - y_k\|^2 - \|u_k - y_k\|^2 - \|z_k - u_k\|^2.$$

Also, for some $c > 0$, we have

$$2r_k \langle -p, x - u_k \rangle = 2r_k \langle -p, x - z_k \rangle + 2r_k \langle -p, z_k - u_k \rangle \\ \le 2r_k \langle -p, x - z_k \rangle + \tfrac{2}{c} r_k^2 \|p\|^2 + \tfrac{c}{2} \|u_k - z_k\|^2.$$

Thus, inequality (7) implies

$$0 \le 2r_k \lambda_k F(z_k,x) - 2\rho r_k \|u_k - x\|^2 + 2r_k \langle -p, x - z_k \rangle \\ + \tfrac{2}{c} r_k^2 \|p\|^2 + \tfrac{c}{2} \|u_k - z_k\|^2 + \|y_k - x\|^2 - \|z_k - x\|^2 \qquad (8) \\ - \|u_k - y_k\|^2 - \|z_k - u_k\|^2 + 2r_k \langle v, x - u_k \rangle.$$

By Lemma 4, for all $k \ge 1$,

$$\begin{aligned} \|y_k - x\|^2 &= \|x_k - \gamma_k(x_k - x_{k-1}) - x\|^2 \\ &= \|(1 - \gamma_k)(x_k - x) + \gamma_k(x_{k-1} - x)\|^2 \qquad (9) \\ &= (1 - \gamma_k)\|x_k - x\|^2 + \gamma_k \|x_{k-1} - x\|^2 - \gamma_k(1 - \gamma_k)\|x_k - x_{k-1}\|^2. \end{aligned}$$

Also, for each $k \ge 1$, choose $b_k \in (\gamma_k, 1)$ such that

$$2\gamma_k \langle u_k - x_k, x_k - x_{k-1} \rangle \ge -\tfrac{\gamma_k}{b_k} \|u_k - x_k\|^2 - \gamma_k b_k \|x_k - x_{k-1}\|^2.$$

Then

$$\begin{aligned} \|u_k - y_k\|^2 &= \|u_k - x_k + \gamma_k(x_k - x_{k-1})\|^2 \\ &= \|u_k - x_k\|^2 + \gamma_k^2 \|x_k - x_{k-1}\|^2 + 2\gamma_k \langle u_k - x_k, x_k - x_{k-1} \rangle \qquad (10) \\ &\ge \left(1 - \tfrac{\gamma_k}{b_k}\right) \|u_k - x_k\|^2 + \gamma_k(\gamma_k - b_k)\|x_k - x_{k-1}\|^2. \end{aligned}$$

Plugging (9) and (10) into (8), and using the monotonicity of $F$, we get

$$0 \leq r_k \lambda_k F(z_k, x) + r_k \lambda_k \left( F(z_k, x) - \tfrac{2}{\lambda_k} \langle p, x - z_k \rangle \right) + \tfrac{1}{c} r_k^2 \|p\|^2$$
$$- (1 + 2\rho r_k)\|z_k - x\|^2 + (1 - \gamma_k)\|x_k - x\|^2 + \gamma_k \|x_{k-1} - x\|^2$$
$$- \gamma_k(1 - b_k)\|x_k - x_{k-1}\|^2 - \left(1 - \tfrac{\gamma_k}{b_k}\right)\|u_k - x_k\|^2$$
$$- (1 - c)\|u_k - z_k\|^2 + 2r_k \langle v, x - u_k \rangle.$$

Since $p \in \mathcal{N}_{\mathbf{S}_F}(x)$, we have

$$\sigma_{\mathbf{S}_F}\left(\tfrac{2}{\lambda_k} p\right) = \left\langle \tfrac{2}{\lambda_k} p, x \right\rangle.$$

Then, by Fenchel–Moreau duality, a term in the first line above rewrites as

$$F(z_k, x) - \tfrac{2}{\lambda_k} \langle p, x - z_k \rangle = \left\langle \tfrac{2}{\lambda_k} p, z_k \right\rangle + F(z_k, x) - \left\langle \tfrac{2}{\lambda_k} p, x \right\rangle$$
$$\leq \sup_{v \in \mathcal{H}} \left\{ \left\langle \tfrac{2}{\lambda_k} p, v \right\rangle + F(v, x) \right\} - \sigma_{\mathbf{S}_F}\left(\tfrac{2}{\lambda_k} p\right)$$
$$= \mathcal{F}_F\left(x, \tfrac{2}{\lambda_k} p\right) - \sigma_{\mathbf{S}_F}\left(\tfrac{2}{\lambda_k} p\right).$$

Combining, we justify (3).

**2.** Using iteration (4i) of Algorithm (RIPSA) and Lemma 4, we deduce

$$\|x_{k+1} - x\|^2 = \|\alpha_k(x_k - x) + \beta_k(z_k - x) + (1 - \alpha_k - \beta_k)(g(x_k) - x)\|^2$$
$$= \alpha_k \|x_k - x\|^2 + \beta_k \|z_k - x\|^2 + (1 - \alpha_k - \beta_k)\|g(x_k) - x\|^2$$
$$- \alpha_k(1 - \alpha_k - \beta_k)\|g(x_k) - x_k\|^2 - \beta_k(1 - \alpha_k - \beta_k)\|g(x_k) - z_k\|^2$$
$$- \alpha_k \beta_k \|z_k - x_k\|^2.$$

Using (3), it follows that

$$\|x_{k+1} - x\|^2 \leq \frac{\beta_k}{1 + 2\rho r_k}\left( \beta_k r_k \lambda_k F(z_k, x) + \beta_k r_k \lambda_k \left[ \mathcal{F}_F\left(x, \tfrac{2}{\lambda_k} p\right) - \sigma_{\mathbf{S}_F}\left(\tfrac{2}{\lambda_k} p\right) \right] \right.$$
$$\left. + \tfrac{2}{c} r_k^2 \beta_k \|p\|^2 \right)$$
$$+ \left( \alpha_k + \tfrac{\beta_k(1 - \gamma_k)}{1 + 2\rho r_k} \right)\|x_k - x\|^2 + \tfrac{\beta_k \gamma_k}{1 + 2\rho r_k}\|x_{k-1} - x\|^2$$
$$- \tfrac{\beta_k}{1 + 2\rho r_k} \gamma_k(1 - b_k)\|x_k - x_{k-1}\|^2 - \tfrac{\beta_k}{1 + 2\rho r_k}\left(1 - \tfrac{\gamma_k}{b_k}\right)\|u_k - x_k\|^2$$
$$- \tfrac{\beta_k}{1 + 2\rho r_k}\left(1 - \tfrac{c}{2}\right)\|u_k - z_k\|^2 + (1 - \alpha_k - \beta_k)\|g(x_k) - x\|^2$$
$$- \alpha_k \beta_k \|z_k - x_k\|^2 + \frac{2r_k \beta_k}{1 + 2\rho r_k} \langle v, x - u_k \rangle,$$

which is exactly the required inequality (4). ∎

## 4 Weak convergence of the algorithm

We first address the weak ergodic convergence of the proximal algorithm (RIPSA) in the case $\rho = \alpha_k = 1 - \beta_k = 0$. These parameters were originally introduced in (RIPSA) to guarantee the strong convergence of the algorithm. In this simplified setting, (RIPSA) reduces to

$$y_k = (1 - \gamma_k)x_k + \gamma_k x_{k-1}, \quad u_k = J_{r_k}^G(y_k), \quad x_{k+1} = J_{r_k}^{\lambda_k F}(u_k). \tag{11}$$

We mean by weak ergodic convergence of the sequence $\{x_k\}$ generated by (RIPSA), the weak convergence of the associate ergodic average $\{\tilde{x}_n\}$ defined by

$$\tilde{x}_n := \frac{1}{\sigma_n} \sum_{k=1}^{n} r_k x_k \text{ where } \sigma_n = \sum_{k=1}^{n} r_k.$$

**Theorem 1** *Suppose the assumptions of Lemma 10 hold, and let $\rho = \alpha_k = 1 - \beta_k = 0$. Assume in addition that condition (H4) is satisfied, that $\{\gamma_k\}$ is nondecreasing, $\{r_k\}$ is nonincreasing, and that*

$$0 < c < 2, \quad \sum_{k=0}^{\infty} r_k = \infty, \quad \sum_{k=0}^{\infty} r_k^2 < \infty. \tag{12}$$

*Then, for the sequence $\{x_k\}$ generated by (RIPSA), cf. (11), we obtain the weak ergodic convergence to a solution $\bar{x}$ of the hierarchical equilibrium problem (HEP). Suppose moreover that*

$$\liminf_{k \to +\infty} \lambda_k r_k > 0,$$

*then the whole sequence $\{x_k\}$ weakly converges to a solution of (HEP).*

**Proof.** To establish weak convergence of $\{x_k\}$, we apply Opial-Passty's Lemma 9.

**Step 1. Convergence of distances.** Let $x \in \mathbf{S}$. By setting $v = 0, \rho = \alpha_k = 1 - \beta_k = 0$ and $\gamma_k < b_k < 1$ in Lemma 10(3), we obtain

$$\begin{aligned}\|x_{k+1} - x\|^2 &\leq r_k \lambda_k F(x_{k+1}, x) + r_k \lambda_k \left[ \mathcal{F}_F\left(x, \tfrac{2}{\lambda_k} p\right) - \sigma_{\mathbf{S}_F}\left(\tfrac{2}{\lambda_k} p\right) \right] \\ &\quad + \tfrac{2}{c} r_k^2 \|p\|^2 + (1 - \gamma_k) \|x_k - x\|^2 + \gamma_k \|x_{k-1} - x\|^2 \\ &\quad - \gamma_k (1 - b_k) \|x_k - x_{k-1}\|^2 - \left(1 - \tfrac{c}{2}\right) \|u_k - x_{k+1}\|^2 \\ &\quad - \left(1 - \tfrac{\gamma_k}{b_k}\right) \|u_k - x_k\|^2. \end{aligned} \tag{13}$$

Setting $a_k = \|x_k - x\|^2$ and using monotonicity of $F$ gives

$$a_{k+1} - (1-\gamma_k)a_k - \gamma_k a_{k-1} + \gamma_k(1-b_k)\|x_k - x_{k-1}\|^2$$
$$+ \left(1 - \tfrac{\gamma_k}{b_k}\right)\|u_k - x_k\|^2 + \left(1 - \tfrac{c}{2}\right)\|u_k - x_{k+1}\|^2 + r_k \lambda_k F(x, x_{k+1})$$
$$\leq r_k \lambda_k \left[\mathcal{F}_F\left(x, \tfrac{2}{\lambda_k}p\right) - \sigma_{\mathbf{S}_F}\left(\tfrac{2}{\lambda_k}p\right)\right] + \tfrac{2}{c} r_k^2 \|p\|^2. \tag{14}$$

Since $1 - b_k > 0$, $b_k - \gamma_k > 0$, and $2 - c > 0$, and by assumption

$$\sum_{k=1}^{\infty} r_k \lambda_k \left[\mathcal{F}_F\left(x, \tfrac{2}{\lambda_k}p\right) - \sigma_{\mathbf{S}_F}\left(\tfrac{2}{\lambda_k}p\right)\right] < \infty, \qquad \sum_{k=1}^{\infty} r_k^2 < \infty,$$

Lemma 7 implies convergence of $\{\|x_k - x\|\}$ as well as of the positive series

$$\sum_{k=1}^{\infty} \gamma_k(1-b_k)\|x_k - x_{k-1}\|^2, \quad \sum_{k=1}^{\infty}\left(1 - \tfrac{\gamma_k}{b_k}\right)\|u_k - x_k\|^2,$$
$$\sum_{k=1}^{\infty} \|u_k - x_{k+1}\|^2, \quad \sum_{k=1}^{\infty} r_k \lambda_k F(x, x_{k+1}).$$

We conclude from $\sum_{k=1}^{\infty} \|u_k - x_{k+1}\|^2 < \infty$, that $\{u_k - x_{k+1}\} \to 0$ strongly in $\mathcal{H}$.

**Step 2. Identification of cluster points for weak ergodic convergence.** Let $x^*$ be a weak cluster point of $\{\tilde{x}_k\}$, with

$$x^* = w-\lim_{n \in I} \tilde{x}_{n+1}$$

along some subsequence $I \subset \mathbb{N}$.

To prove that $x^*$ is a solution of the hierarchical problem (**HEP**), it suffices to show that

$$0 \in A^G(x^*) + \mathcal{N}_{\mathbf{S}_F}(x^*). \tag{15}$$

Since $\mathbb{R}_+(C - \mathbf{S}_F)$ is a closed linear subspace of $\mathcal{H}$, the operator $A^G + \mathcal{N}_{\mathbf{S}_F}$ is maximally monotone (see [37, 38]). Thus, proving (15) reduces to showing that

$$\langle v, x - x^* \rangle \geq 0, \quad \text{for all } (x, v) \in \text{gra}(A^G + \mathcal{N}_{\mathbf{S}_F}). \tag{16}$$

Take $x \in \mathbf{S}_F$ and $p \in \mathcal{N}_{\mathbf{S}_F}(x)$ such that $v - p \in A^G(x)$. Setting $\rho = \alpha_k = 1 - \beta_k = 0$ and $\gamma_k < b_k < 1$, inequality (4) in Lemma 10 yields

$$\|x_{k+1} - x\|^2 \leq r_k \lambda_k F(z_k, x) + r_k \lambda_k \left[\mathcal{F}_F\left(x, \tfrac{2}{\lambda_k}p\right) - \sigma_{\mathbf{S}_F}\left(\tfrac{2}{\lambda_k}p\right)\right]$$
$$+ \tfrac{2}{c} r_k^2 \|p\|^2 + (1-\gamma_k)\|x_k - x\|^2 + \gamma_k \|x_{k-1} - x\|^2$$
$$- \gamma_k(1-b_k)\|x_k - x_{k-1}\|^2 - \left(1 - \tfrac{\gamma_k}{b_k}\right)\|u_k - x_k\|^2$$
$$- \left(1 - \tfrac{c}{2}\right)\|u_k - z_k\|^2 + 2r_k \langle v, x - u_k \rangle.$$

Setting $a_k = \|x_k - x\|^2$, discarding nonpositive terms, and assuming that $\{\gamma_k\}$ is nonincreasing, we obtain

$$2r_k \langle v, u_k - x \rangle \leq (a_k - a_{k+1}) + \gamma_{k-1} a_{k-1} - \gamma_k a_k \\ + r_k \lambda_k \left[ \mathcal{F}_F\!\left(x, \tfrac{2}{\lambda_k} p\right) - \sigma_{\mathbf{S}_F}\!\left(\tfrac{2}{\lambda_k} p\right) \right] + \tfrac{2}{c} \|p\|^2 r_k^2.$$

Summing for $k = 1, \ldots, n$ and dividing by $\sigma_n = \sum_{k=1}^n r_k$, we deduce

$$\langle v, \tilde{u}_n - x \rangle = \frac{1}{\sigma_n} \sum_{k=1}^n r_k \langle v, u_k - x \rangle \leq \frac{1}{2\sigma_n} \left( a_1 + \gamma_0 a_0 + \bar{M} \right), \tag{17}$$

where

$$\bar{M} = \tfrac{1}{c} \|p\|^2 \sum_{k=1}^\infty r_k^2 + \sum_{k=1}^\infty r_k \lambda_k \left[ \mathcal{F}_F\!\left(x, \tfrac{2}{\lambda_k} p\right) - \sigma_{\mathbf{S}_F}\!\left(\tfrac{2}{\lambda_k} p\right) \right] < +\infty.$$

Since $\sigma_n \to +\infty$ as $n \to +\infty$, it follows that

$$\limsup_{n \to +\infty} \langle v, \tilde{u}_n - x \rangle \leq 0.$$

Moreover, using $x^* = w\text{-}\lim_{n \in I} \tilde{x}_{n+1}$ and the fact that $u_k - x_{k+1} \to 0$ strongly, we deduce that $x^* = w\text{-}\lim_{n \in I} \tilde{u}_n$. Hence

$$\langle v, x^* - x \rangle = \lim_{n \in I, n \to \infty} \langle v, \tilde{u}_n - x \rangle \leq \limsup_{n \to \infty} \langle v, \tilde{u}_n - x \rangle \leq 0,$$

which proves (16). Thus, $x^*$ is a solution of the hierarchical problem (**HEP**).

By Lemma 9, the entire sequence $\{\tilde{x}_n\}$ converges weakly to a solution of (**HEP**). This establishes the weak ergodic convergence result.

**Step 3. Identification of cluster points for weak convergence.** Let $x^*$ be a weak cluster point of $\{x_k\}$, and consider a subsequence $\{x_{k+1} : k \in I \subset \mathbb{N}\}$ such that

$$x^* = w\text{-}\lim_{k \in I} x_{k+1}.$$

Step 3.1: Proof that $x^* \in \mathbf{S}_F$. Fix $y \in C$. From the iteration (11), we have $x_{k+1} = J_{r_k}^{\lambda_k F}(u_k)$, which yields

$$F(x_{k+1}, y) + \tfrac{1}{r_k \lambda_k} \langle x_{k+1} - u_k, y - x_{k+1} \rangle \geq 0.$$

Since $F$ is monotone and weakly lower semicontinuous in its second argument, it follows that

$$F(y,x^*) \le \liminf_{k\in I, k\to\infty} F(y,x_{k+1}) \le \liminf_{k\in I, k\to\infty} -F(x_{k+1},y)$$
$$\le \limsup_{k\in I, k\to\infty} \frac{1}{r_k\lambda_k}\|x_{k+1}-u_k\|\cdot\|y-x_{k+1}\| = 0.$$

The last equality holds since $\liminf_{k\to\infty} r_k\lambda_k > 0$, $\{x_{k+1} : k \in I\}$ is bounded, and $\|u_k - x_{k+1}\| \to 0$. By Lemma 1, we conclude that $x^* \in \mathbf{S}_F$.

<u>Step 3.2:</u> We prove that $x^* \in \mathbf{S}$. Recall that the iterations are given by $x_{k+1} = J_{r_k}^{\lambda_k F}(u_k)$, $u_k = J_{r_k}^G(y_k)$. By monotonicity of $F$ and $G$, we obtain, for all $x \in \mathbf{S}_F$,

$$0 \le r_k\lambda_k F(x,x_{k+1}) \le -r_k\lambda_k F(x_{k+1},x) \le \langle x_{k+1}-u_k, x-x_{k+1}\rangle$$
$$= \|u_k-x\|^2 - \|u_k-x_{k+1}\|^2 - \|x_{k+1}-x\|^2, \tag{18}$$

$$2r_k G(x,u_k) \le -2r_k G(u_k,x) \le 2\langle u_k-y_k, x-u_k\rangle$$
$$= \|y_k-x\|^2 - \|u_k-y_k\|^2 - \|u_k-x\|^2. \tag{19}$$

From Lemma 4 and (18), we deduce

$$\|x_{k+1}-x\|^2 \le \|u_k-x\|^2 - \|x_{k+1}-u_k\|^2,$$

together with

$$\|y_k-x\|^2 = (1-\gamma_k)\|x_k-x\|^2 + \gamma_k\|x_{k-1}-x\|^2 - \gamma_k(1-\gamma_k)\|x_k-x_{k-1}\|^2.$$

Moreover, as in (10), for $b > 0$ we obtain

$$\|u_k - y_k\|^2 \ge (\gamma_k - b)\left(\gamma_k\|x_k - x_{k-1}\|^2 - \tfrac{1}{b}\|u_k - x_k\|^2\right).$$

Substituting these relations into (19), we arrive at

$$2r_k G(x,u_k) \le (1-\gamma_k)\|x_k-x\|^2 + \gamma_k\|x_{k-1}-x\|^2 - \gamma_k(1-\gamma_k)\|x_k-x_{k-1}\|^2$$
$$- \gamma_k(\gamma_k-b)\|x_k-x_{k-1}\|^2 - \left(1-\tfrac{\gamma_k}{b}\right)\|u_k-x_k\|^2$$
$$- \|x_{k+1}-x\|^2 - \|x_{k+1}-u_k\|^2$$
$$= \left(\|x_k-x\|^2 - \|x_{k+1}-x\|^2\right) + \gamma_k\left(\|x_{k-1}-x\|^2 - \|x_k-x\|^2\right)$$
$$- \left(\gamma_k(1-b)\|x_k-x_{k-1}\|^2 + \|u_k-x_{k+1}\|^2 + \left(1-\tfrac{\gamma_k}{b}\right)\|u_k-x_k\|^2\right).$$

Since $\gamma_k \le b_k$, this gives

$$2r_k G(x,u_k) \le (a_k - a_{k+1}) + \gamma_k(a_{k-1} - a_k),$$

where $a_k = \|x_k - x\|^2$. Assuming $\{\gamma_k\}$ is nonincreasing, we further obtain

$$2r_k G(x,u_k) + b_k \le (a_k - a_{k+1}) + \gamma_{k-1}a_{k-1} - \gamma_k a_k.$$

Summing over $k \ge 1$, we deduce

$$2\sum_{k=1}^{\infty} r_k G(x, u_k) + \sum_{k=1}^{\infty} b_k \leq a_1 + \gamma_0 a_0 < +\infty.$$

By Lemma 5, since $\sum_{k=1}^{\infty} r_k = +\infty$, it follows that

$$\liminf_{k \to \infty} G(x, u_k) \leq 0.$$

Finally, since $\{x_k\}$ converges weakly to $x^*$ and $G$ is weakly lower semicontinuous in its second argument, we obtain

$$G(x, x^*) \leq \liminf_{k \to \infty} G(x, u_k) \leq 0.$$

This holds for all $x \in \mathbf{S}_F$.

Therefore, $x^*$ is a solution to the hierarchical problem (**HEP**). By Lemma 9, we conclude that the entire sequence $\{x_k\}$ converges weakly to a solution $\bar{x}$ of (**HEP**). ∎

## 5 Strong convergence of the algorithm

To analyze the strong convergence of the algorithm (RIPSA), we proceed in two different ways:

- The first approach strengthens the monotonicity assumption: using the strong monotonicity of the bifunction $G$, we establish the strong convergence of the sequence $\{x_k\}$ to the unique solution of problem (**EP**).
- The second approach perturbs the equilibrium problem by introducing a contraction mapping $g$, ensuring that the algorithm is executable, and subsequently proving the strong convergence of $\{x_k\}$ to a solution of problem (**EP**), with a specific selection of this solution relative to $g$.

### 5.1 Strong convergence for strongly monotone bifunctions

As a distinctive feature, the set of solutions $\mathbf{S}_F$ of a strongly monotone equilibrium problem is nonempty and reduces to a singleton $\{\bar{x}\}$. We now prove the strong convergence of the sequence $\{x_k\}$ generated by (RIPSA) when $G$ is $\rho$-strongly monotone and $\alpha_k = 1 - \beta_k = 0$. In this case, the solution $\bar{x}$ of (**HEP**) is unique.

**Theorem 2** *Suppose that the assumptions of Lemma 10 hold with $\alpha_k = 1 - \beta_k = 0$, together with condition* (*H*5). *Assume further that $G$ is $\rho$-strongly monotone, with $\rho > 0$. Suppose moreover that*

$$\gamma_k \in [0, 1/2], \quad \limsup_{k \to \infty} \lambda_k \left[ \mathcal{F}_F\left(\bar{x}, \tfrac{2}{\lambda_k} p\right) - \sigma_{\mathbf{S}_F}\left(\tfrac{2}{\lambda_k} p\right) \right] \leq 0,$$

$$\sum_k r_k = +\infty, \quad \sum_k r_k^2 < +\infty, \quad \gamma_k < b_k < 1, \quad 0 < c < 2. \tag{20}$$

*Then, the solution set $S_F$ of* (**EP**) *is nonempty, the solution set $S$ of* (**HEP**) *is a singleton $\{\bar{x}\}$, and the sequence $\{x_k\}$ generated by (RIPSA) converges strongly to $\bar{x}$.*

**Proof.** Existence and uniqueness follow from [9, Thm. 4.3] and the strong monotonicity of $G$. Let $\bar{x}$ denote the unique element of **S**.

To prove strong convergence of $\{x_k\}$ to $\bar{x}$, we return to Lemma 10. Using $\rho > 0$, $\alpha_k = 1 - \beta_k = 0$, $\gamma_k < b_k < 1$, and $c < 2$, and noting that

$$\frac{1}{1 + 2\rho r_k} = 1 - \frac{2\rho r_k}{1 + 2\rho r_k},$$

we deduce from (4) that

$$\begin{aligned}
\|x_{k+1} - \bar{x}\|^2 &\leq \frac{r_k \lambda_k}{1 + 2\rho r_k} \left[ \mathcal{F}_F\left(\bar{x}, \tfrac{2}{\lambda_k} p\right) - \sigma_{\mathbf{S}_F}\left(\tfrac{2}{\lambda_k} p\right) \right] \\
&\quad + \frac{2 r_k^2}{c} \|p\|^2 + \left( 1 - \frac{2\rho r_k}{1 + 2\rho r_k} - \frac{\gamma_k}{1 + 2\rho r_k} \right) \|x_k - \bar{x}\|^2 \\
&\quad + \frac{\gamma_k}{1 + 2\rho r_k} \|x_{k-1} - \bar{x}\|^2.
\end{aligned} \tag{21}$$

In Lemma 8, set

$$a_k = \|x_k - \bar{x}\|^2, \quad t_k = \frac{2\rho r_k}{1 + 2\rho r_k}, \quad \delta_k = \frac{\gamma_k}{1 + 2\rho r_k}, \quad \theta_k = \frac{2\|p\|^2}{c} r_k^2,$$

$$s_k = \frac{\lambda_k}{2\rho} \left[ \mathcal{F}_F\left(\bar{x}, \tfrac{2}{\lambda_k} p\right) - \sigma_{\mathbf{S}_F}\left(\tfrac{2}{\lambda_k} p\right) \right].$$

The first three conditions in (2) follow directly from (20). The last condition, $\sum_k t_k = +\infty$, follows from $\sum_k r_k^2 < +\infty$ and $\sum_k r_k = +\infty$, since $r_k \to 0$ implies $r_k \leq \tfrac{1}{2}$ for large $k$, and hence

$$\sum_k t_k = \sum_k \frac{2\rho r_k}{1 + 2\rho r_k} \geq \frac{2\rho}{1 + \rho} \sum_k r_k = +\infty.$$

Lemma 8 therefore ensures that $a_k \to 0$, i.e., $\{x_k\}$ converges strongly to $\bar{x}$. ∎

## 5.2 Strong convergence for non-strongly monotone bifunctions

**Theorem 3** *Suppose the assumptions of Lemma 10 and condition* (H4) *hold, and let $g$ be a $\delta$-contraction on $C$, with $0 \leq \delta < 1$. Assume moreover that*

$$\lim_{k\to\infty} \alpha_k < 1, \quad \lim_{k\to\infty} \beta_k = 0, \quad \lim_{k\to\infty} \lambda_k = +\infty, \quad \lim_{k\to\infty} \mu_k = 0,$$

$$\lim_{k\to\infty} \frac{\beta_k}{\mu_k}\omega_k = 0, \quad \lim_{k\to\infty} r_k > 0, \quad \sum_{k=0}^{\infty} \beta_k \omega_k < +\infty, \quad \sum_{k=0}^{\infty} \mu_k = +\infty. \quad (22)$$

*Then the sequence $\{x_k\}$ generated by (RIPSA) converges strongly to a selected solution $\bar{x}$ of (**HEP**); more precisely, $\bar{x} \in S$ is the unique solution of the variational inequality*

$$\langle \bar{x} - G(\bar{x}), \bar{x} - x \rangle \leq 0, \quad \text{for all } x \in S. \quad (23)$$

**Proof.** The proof is divided into several steps, where the title of each step indicates its main purpose.

**Step 1: Existence and uniqueness of a solution to** (23)**.** Inequality (23) is precisely the variational characterization of the metric projection

$$\bar{x} = \Pi_S(G(\bar{x})),$$

where $\Pi_S(G(\bar{x}))$ denotes the unique point in $S$ closest (in the metric sense) to $G(\bar{x})$. Since $\Pi_{S_F}$ is nonexpansive and $g$ is a contraction on $C$, their composition

$$\Pi_{S_F} \circ g$$

is a contraction. By Banach's fixed-point theorem, this contraction admits a unique fixed point $\bar{x}$.

**Step 2: Boundedness of the sequences** $\{x_k\}, \{u_k\}, \{z_k\}$**.**

Boundedness of $\{x_k\}$. We first establish boundedness of the real sequence $\{\|x_k - \bar{x}\|^2\}$.

Define

$$\omega_k := \lambda_k \left[ \mathcal{F}_F\left(\bar{x}, \tfrac{2}{\lambda_k} p\right) - \sigma_{S_F}\left(\tfrac{2}{\lambda_k} p\right) \right] + \frac{2}{c} r_k \|p\|^2, \quad (24)$$

and set

$$M := \max\left\{ \|x_0 - \bar{x}\|^2, \|x_1 - \bar{x}\|^2, \frac{2\|g(\bar{x}) - \bar{x}\|^2}{(1-\delta)^2} \right\},$$

where $x_0, x_1$ are the initial points of (RIPSA), and $0 \leq \delta < 1$ is the contraction constant of $g$.

Suppose $\|x_j - \bar{x}\|^2 \leq M$ for all $0 \leq j \leq k$. We show by induction that $\|x_{k+1} - \bar{x}\|^2 \leq M$.

Let $\mu_k = 1 - \alpha_k - \beta_k$. Rewriting inequality (4) in Lemma 10 with $\gamma_k < b_k < 1$ and $\rho = 0$ (since $G$ is only assumed monotone), we obtain

$$\|x_{k+1} - \bar{x}\|^2 \leq \beta_k r_k \lambda_k F(z_k, \bar{x}) + \beta_k r_k \lambda_k \left[ \mathcal{F}_F\left(x, \tfrac{2}{\lambda_k} p\right) - \sigma_{\mathbf{S}_F}\left(\tfrac{2}{\lambda_k} p\right) \right]$$
$$+ \tfrac{2}{c} r_k^2 \beta_k \|p\|^2 + (\alpha_k + \beta_k(1-\gamma_k)) \|x_k - \bar{x}\|^2 + \beta_k \gamma_k \|x_{k-1} - \bar{x}\|^2$$
$$- \beta_k \gamma_k (1 - b_k) \|x_k - x_{k-1}\|^2 - \beta_k \left(1 - \tfrac{c}{2}\right) \|u_k - z_k\|^2 \quad (25)$$
$$- \beta_k \left(1 - \tfrac{\gamma_k}{b_k}\right) \|u_k - x_k\|^2 + \mu_k \|g(x_k) - \bar{x}\|^2 - \alpha_k \beta_k \|z_k - x_k\|^2.$$

Using the contraction property of $g$,

$$\|g(x_k) - \bar{x}\|^2 \leq \delta \|x_k - \bar{x}\|^2 + \frac{1}{1-\delta} \|g(\bar{x}) - \bar{x}\|^2.$$

Since $\gamma_k < b_k < 1$ and $0 < c < 2$, inequality (25) yields

$$\|x_{k+1} - \bar{x}\|^2 \leq r_k \beta_k \left( \lambda_k \left[ \mathcal{F}_F\left(x, \tfrac{2}{\lambda_k} p\right) - \sigma_{\mathbf{S}_F}\left(\tfrac{2}{\lambda_k} p\right) \right] + \tfrac{2}{c} r_k \|p\|^2 \right)$$
$$+ (\alpha_k + \beta_k(1-\gamma_k) + \mu_k \delta) \|x_k - \bar{x}\|^2 + \beta_k \gamma_k \|x_{k-1} - \bar{x}\|^2$$
$$+ \frac{\mu_k}{1-\delta} \|g(\bar{x}) - \bar{x}\|^2.$$

By the induction hypothesis and using $\alpha_k + \beta_k + \mu_k = 1$, it follows that

$$\|x_{k+1} - \bar{x}\|^2 \leq r_k \beta_k \omega_k + \left( \alpha_k + \beta_k(1-\gamma_k) + \mu_k \delta + \beta_k \gamma_k \right) M + \mu_k (1-\delta)^2$$
$$= M + \left[ \frac{r_k \beta_k \omega_k}{\mu_k} - (1-\delta)\left(M - \frac{1}{(1-\delta)^2} \|g(\bar{x}) - \bar{x}\|^2\right) \right] \mu_k$$
$$\leq M.$$

The last inequality holds because $M \geq \frac{1}{(1-\delta)^2} \|g(\bar{x}) - \bar{x}\|^2$ and

$$\lim_{k \to \infty} \frac{r_k \beta_k \omega_k}{\mu_k} = 0,$$

which follows from the assumptions $\lim_{k \to \infty} \frac{\beta_k}{\mu_k} \omega_k = 0$ and $\lim_{k \to \infty} r_k < 1$.

Thus, by induction, $\{\|x_k - \bar{x}\|^2\}$ is bounded by $M$, hence $\{x_k\}$ is bounded.

<u>Boundedness of $\{u_k\}$.</u> Let $\bar{y}_k := J_{r_k}^G(\bar{x}_0)$ for some $\bar{x}_0 \in C$. By Lemma 3(i),

$$\|\bar{y}_k\| \leq \|\bar{x}_0\| + r_k \|A_o^G G(\bar{x}_0)\|.$$

Then, by Lemma 3(iii),

$$\|u_k\| \leq (1-\gamma_k) \|x_k - \bar{x}_0\| + \gamma_k \|x_{k-1} - \bar{x}_0\| + \|\bar{x}_0\| + r_k \|A_o^G(\bar{x}_0)\|.$$

Since the right-hand side is bounded, $\{u_k\}$ is bounded.

<u>Boundedness of $\{z_k\}$.</u> Since $z_k = J_{r_k}^{\lambda_k F}(u_k)$ and $\bar{x} \in \mathbf{S}_F$, we have $\bar{x} = J_{r_k}^{\lambda_k F}(\bar{x})$. By Lemma 3(iii),

$$\|z_k - \bar{x}\| \leq \|u_k - \bar{x}\|,$$

so boundedness of $\{u_k\}$ implies boundedness of $\{z_k\}$.

**Step 3: Strong convergence of $\{x_k\}$ towards $\bar{x}$.**

Applying $\|a+b\|^2 \leq \|a\|^2 + 2\langle b, a+b\rangle$ and convexity of $\|\cdot\|^2$, we deduce

$$\begin{aligned}\|x_{k+1} - \bar{x}\|^2 &= \|\alpha_k(x_k - \bar{x}) + \beta_k(z_k - \bar{x}) + \mu_k(g(x_k) - \bar{x})\|^2 \\ &\leq \|\alpha_k(x_k - \bar{x}) + \beta_k(z_k - \bar{x})\|^2 + 2\mu_k\langle g(x_k) - \bar{x}, x_{k+1} - \bar{x}\rangle \\ &\leq (\alpha_k + \beta_k)^2 \left(\tfrac{\alpha_k}{\alpha_k+\beta_k}\|x_k - \bar{x}\|^2 + \tfrac{\beta_k}{\alpha_k+\beta_k}\|z_k - \bar{x}\|^2\right) \\ &\quad + 2\mu_k\langle g(x_k) - g(\bar{x}), x_{k+1} - \bar{x}\rangle + 2\mu_k\langle g(\bar{x}) - \bar{x}, x_{k+1} - \bar{x}\rangle.\end{aligned}$$

Using (4) and that $g$ is a $\delta$-contraction, we obtain

$$\begin{aligned}\|x_{k+1} - \bar{x}\|^2 &\leq (\alpha_k + \beta_k)\bigl(\alpha_k\|x_k - \bar{x}\|^2 + \beta_k\|z_k - \bar{x}\|^2\bigr) \\ &\quad + \delta^2\mu_k\|x_k - \bar{x}\|^2 + \mu_k\|x_{k+1} - \bar{x}\|^2 + 2\mu_k\langle g(\bar{x}) - \bar{x}, x_{k+1} - \bar{x}\rangle \\ &\quad + r_k(\alpha_k + \beta_k)\beta_k\Bigl(\lambda_k\bigl[\mathcal{F}_F\bigl(x, \tfrac{2}{\lambda_k}p\bigr) - \sigma_{\mathbf{S}_F}\bigl(\tfrac{2}{\lambda_k}p\bigr)\bigr] + \tfrac{2}{c}r_k\|p\|^2\Bigr).\end{aligned}$$

Setting $a_k := \|x_k - \bar{x}\|^2$, $g_k := \langle g(\bar{x}) - \bar{x}, x_{k+1} - \bar{x}\rangle$, and recalling (24), we arrive at

$$a_{k+1} \leq \bigl(\alpha_k - \beta_k\gamma_k - \mu_k(1 - \tfrac{\delta^2}{1-\mu_k})\bigr)a_k + \beta_k\gamma_k a_{k-1} + \tfrac{2\mu_k}{1-\mu_k}g_k + r_k\beta_k\omega_k. \qquad (26)$$

By setting

$$\delta_k = \beta_k\gamma_k, \quad \theta_k = r_k\beta_k\omega_k, \quad t_k = \mu_k\Bigl(1 - \tfrac{\delta^2}{1-\mu_k}\Bigr), \quad s_k = \tfrac{2\mu_k g_k}{1-\mu_k},$$

we verify that the assumptions of Lemma 8 are satisfied:

- Since $\beta_k \to 0$ and $\gamma_k \leq 1$, we have $\delta_k \to 0$, hence eventually $\delta_k \in [0, \tfrac{1}{2}]$.
- By boundedness of $\{r_k\}$ and assumption $\sum \beta_k\omega_k < \infty$, it follows that $\sum \theta_k < \infty$.
- Because $\sum \mu_k = \infty$ and $\mu_k \to 0$, the term $1 - \tfrac{\delta^2}{1-\mu_k}$ stays bounded below by a positive constant. Hence $\sum t_k = +\infty$.
- Since $\mu_k \to 0$, we have $\tfrac{2}{1-\mu_k} \leq 4$ for $k$ large. Boundedness of $\{x_k\}$ implies $|g_k| \leq A_0$, for some $A_0 > 0$, so $\limsup_{k\to\infty} s_k \leq 0$.

All conditions of Lemma 8 are satisfied. Therefore, $a_k \to 0$, i.e.,

$$x_k \to \bar{x} = \Pi_{\mathbf{S}}(G(\bar{x})) \quad \text{strongly.}$$

∎

# 6 Comments on Particular Cases

The aim of this section is to sketch and discuss some particular cases of equilibrium problems, such as optimization, variational inequalities, and fixed-point problems.

## 6.1 Minimization Problems

Let $\psi : \mathcal{H} \to \mathbb{R} \cup \{+\infty\}$ be a convex, lower semicontinuous function whose effective domain $\mathrm{dom}(\psi)$ contains the convex closed subset $C$. We consider the minimization problem

$$\min_{x \in C} \psi(x).$$

There are two classical ways to represent this minimization problem as an equilibrium problem:

**First way.**

The first approach is to define one of the bifunctions, say $F(x,y)$ or $G(x,y)$, as

$$G(x,y) := \psi(y) - \psi(x).$$

In this case, assumptions $(H1)$ and $(C5)$ are satisfied, and the associated bifunction is monotone. Therefore, if we take $G(x,y) = \psi(y) - \psi(x)$, Theorem 1 or Theorem 3 can be applied to deduce, respectively, weak or strong convergence to a solution of the hierarchical problem (**HEP**).

However, this choice is not suitable for Theorem 2, since the strong monotonicity of $G$ cannot be verified. Thus, an alternative formulation of $G$ is recommended (see Second way below).

Now, for $G(x,y) := \psi(y) - \psi(x)$, let us compute the form of $J_r^G(x)$ for all $x \in \mathcal{H}$. From Lemma 2, we have

$$\begin{aligned}
z = J_r^G(x) &\iff rG(z,y) + \langle y-z, z-x \rangle \geq 0, \text{ for all } y \in C \\
&\iff (\psi + \delta_C)(y) + \left\langle \tfrac{1}{r}(z-x), y \right\rangle \geq (\psi + \delta_C)(z) + \left\langle \tfrac{1}{r}(z-x), z \right\rangle, \text{ for all } y \in \mathcal{H} \\
&\iff 0 \in \partial(\psi + \delta_C)(z) + \tfrac{1}{r}(z-x) \\
&\iff x \in z + r\partial(\psi + \delta_C)(z) \\
&\iff z = \mathrm{prox}_{r(\psi + \delta_C)}(x).
\end{aligned}$$

Next, let us discuss condition $(H4)$ for the special case where the lower-level problem is also a minimization, with

$$F(x,y) = \varphi(y) - \varphi(x), \qquad \varphi : C \to \mathbb{R} \text{ convex l.s.c.}$$

Let $\bar{\varphi} : \mathcal{H} \to \mathbb{R} \cup \{+\infty\}$ denote the extension of $\varphi$ to $\mathcal{H}$ by setting it equal to $+\infty$ outside $C$. Then

$$\mathbf{S}_F = \operatorname*{argmin}_C \varphi = \operatorname*{argmin}_{\mathcal{H}} \bar{\varphi},$$

and $F$ clearly satisfies all conditions of $(H1)$.

For simplicity, assume $\min_C \varphi = 0$. Then

$$\bar{\varphi}(x) - \delta_{\operatorname{argmin}_{\mathcal{H}} \bar{\varphi}}(x) \leq 0, \qquad \text{for all } x \in \mathcal{H}.$$

By Fenchel conjugacy, this implies

$$\bar{\varphi}^*(p) - \sigma_{\operatorname{argmin}_{\mathcal{H}} \bar{\varphi}}(p) \geq 0, \qquad \text{for all } p \in \mathcal{H}.$$

In this case, condition $(H4)$ reduces to the condition $(H_4)$ proposed in [39], namely: for all $u \in \operatorname{argmin}_{\mathcal{H}} \bar{\varphi}$, $p \in \mathcal{N}_{\operatorname{argmin}_{\mathcal{H}} \bar{\varphi}}(u)$,

$$0 \leq \sum_k r_k \lambda_k \left[ \bar{\varphi}^* \left( \frac{2p}{\lambda_k} \right) - \sigma_{\operatorname{argmin}_{\mathcal{H}} \bar{\varphi}} \left( \frac{2p}{\lambda_k} \right) \right] < +\infty. \qquad (27)$$

As an example, if $\varphi(x) = \frac{1}{2} d(x, K)^2$, where $K \subset C$ is a nonempty closed convex set, then for all $p \in \mathcal{H}$,

$$\bar{\varphi}^*(p) - \sigma_{\operatorname{argmin}_{\mathcal{H}} \bar{\varphi}}(p) = \tfrac{1}{2} \|p\|^2.$$

Hence, $(H4)$ (resp. $(H5)$) is equivalent to

$$\sum_k \frac{r_k}{\lambda_k} < +\infty \qquad \left( \text{resp. } \lim_k \frac{r_k}{\lambda_k} = 0 \right).$$

More generally, if $\varphi(x) \geq \frac{1}{2} d(x, K)^2$, then for all $p \in \mathcal{R}(\mathcal{N}_{\operatorname{argmin}_{\mathcal{H}} \bar{\varphi}})$,

$$\bar{\varphi}^*(p) - \sigma_{\operatorname{argmin}_{\mathcal{H}} \bar{\varphi}}(p) \leq \tfrac{1}{2} \|p\|^2,$$

so $(H4)$ still holds whenever $\sum_k \frac{r_k}{\lambda_k} < +\infty$.

**Second way.**

The second approach is to define $G$ by the directional derivative of $\psi$ at $x$ in the direction $y - x$, i.e.,

$$G(x, y) := \psi'(x; y - x) := \lim_{t \to 0^+} \frac{1}{t} \big( \psi(x + t(y - x)) - \psi(x) \big).$$

When $\psi : \mathcal{H} \to \mathbb{R} \cup \{+\infty\}$ is convex, l.s.c., and $\operatorname{dom}(\psi)$ has nonempty interior containing $C$, then for every $(x, y) \in C \times C$,

$$G(x,y) := \psi'(x; y-x) = \sup_{\xi \in \partial \psi(x)} \langle \xi, y-x \rangle,$$

where $\partial \psi$ is the convex subdifferential of $\psi$.

In this case, all conditions $(C1)$–$(C5)$ are satisfied by $G$. Moreover, since $C$ is contained in the interior of $\mathrm{dom}(\psi)$, any solution $\bar{x}$ of (**HEP**) satisfies that $\psi$ is continuous at $\bar{x}$, and $\partial \psi(\bar{x})$ is nonempty, convex, and bounded. Hence, there exists $\xi^* \in \partial \psi(\bar{x})$ such that

$$\langle \xi^*, y - \bar{x} \rangle \geq 0, \quad \text{for all } y \in \mathbf{S}_F.$$

Equivalently,

$$0 \in \xi^* + \partial(\delta_{\mathbf{S}_F})(\bar{x}) \subset \partial \psi(\bar{x}) + \partial(\delta_{\mathbf{S}_F})(\bar{x}) = \partial(\psi + \delta_{\mathbf{S}_F})(\bar{x}),$$

which means $\bar{x}$ minimizes $\psi$ over $\mathbf{S}_F$. Thus $\bar{x}$ is a solution of the hierarchical problem $\min_{\mathbf{S}_F} \psi$.

For this choice of $G$, both Theorems 1 and 3 apply. If, moreover, $\psi$ is strongly convex on $C$, i.e., for some $m > 0$, it holds

$$\psi(tx + (1-t)y) \leq t\psi(x) + (1-t)\psi(y) - mt(1-t)\|x - y\|^2, \text{ for all } x, y \in C, \text{ for all } t \in [0,1],$$

then $G$ is strongly monotone, so Theorem 2 also applies. This was not possible in First way with $G(x,y) = \psi(y) - \psi(x)$.

Finally, note that when

$$G(x,y) = \psi'(x; y-x) = \max_{\xi \in \partial \psi(x)} \langle \xi, y-x \rangle,$$

we obtain

$$\begin{aligned}
z = J_r^G(x) &\iff r \max_{\xi \in \partial \psi(z)} \langle \xi, y-z \rangle + \langle z-x, y-z \rangle \geq 0, \quad \text{for all } y \in C \\
&\iff \exists \xi^* \in \partial \psi(z) \text{ s.t. } \langle r\xi^* + z - x, y - z \rangle \geq 0, \quad \text{for all } y \in C \\
&\iff \exists \xi^* \in \partial \psi(z) \text{ s.t. } -(r\xi^* + z - x) \in \mathcal{N}_C(z) \\
&\iff 0 \in r\partial(\psi + \delta_C)(z) + (z - x) \\
&\iff z = (I + r\partial(\psi + \delta_C))^{-1}(x) = \mathrm{prox}_{r(\psi + \delta_C)}(x).
\end{aligned}$$

In the second equivalence we apply [40, Lemma 3.1], since $\partial \psi(z)$ is convex and weakly compact, $C$ is convex, and the function $p(y, \xi) = \langle r\xi + z - x, y - z \rangle$ is convex in $y$ and concave, upper semicontinuous in $\xi$. The fourth equivalence follows from the identity

$$\partial(\psi + \delta_C)(z) = \partial \psi(z) + \partial(\delta_C)(z) = \partial \psi(z) + \mathcal{N}_C(z),$$

since $\psi$ is continuous on $C$.

## 6.2 Variational inequalities and fixed point problems

### A. Variational inequalities

Let $A : \mathcal{H} \rightrightarrows \mathcal{H}$ be a set-valued mapping such that $Ax$ is nonempty, convex, and weakly compact for each $x \in \mathcal{H}$. We seek to find

$$\bar{x} \in \mathcal{H}, \ \xi^* \in A(\bar{x}) \quad \text{such that} \quad \langle \xi^*, y - \bar{x} \rangle \geq 0 \quad \text{for all } y \in \mathcal{H}. \tag{28}$$

Defining
$$F(x, y) := \max_{\xi \in A(x)} \langle \xi, y - x \rangle,$$

we see that $\bar{x}$ solves (28) if and only if it is a solution of the associated equilibrium problem (**EP**). In this case, $F$ satisfies conditions $(C1)$–$(C3)$ and $(C5)$. Moreover, the $\rho$-monotonicity of $G$ is equivalent to that of the operator $A$:

$$\langle \xi - \eta, x - y \rangle \geq \rho \|x - y\|^2, \quad \text{for all } x, y \in \mathcal{H}, \ \xi \in A(x), \ \eta \in A(y).$$

As in Second way of the previous subsection, we obtain:

$$z = J_r^G(x) \iff r \max_{\xi \in A(z)} \langle \xi, y - z \rangle + \langle z - x, y - z \rangle \geq 0, \text{ for all } y \in \mathcal{H}$$
$$\iff \exists \xi^* \in \partial \psi(z), \ \langle r\xi^* + z - x, y - z \rangle \geq 0, \text{ for all } y \in \mathcal{H}$$
$$\iff x \in z + rA(z)$$
$$\iff z = (I + rA)^{-1}(x) = J_r^A(x).$$

### B. Fixed point problems

Now consider a set-valued mapping $T : \mathcal{H} \rightrightarrows \mathcal{H}$ with nonempty, convex, and weakly compact values. We are interested in the fixed point problem:

$$\text{find} \quad \bar{x} \in \mathcal{H} \quad \text{such that} \quad \bar{x} \in T(\bar{x}). \tag{29}$$

Setting
$$F(x, y) := \max_{\xi \in T(x)} \langle x - \xi, y - x \rangle,$$

we see that $\bar{x}$ solves (**EP**) if and only if $\bar{x}$ solves (29).

In this case, $F$ satisfies conditions $(C1)$–$(C3)$, while condition $(C4)$ is equivalent to

$$\langle \xi - \eta, x - y \rangle \leq \|x - y\|^2, \quad \text{for all } x, y \in C, \ \xi \in T(x), \ \eta \in T(y),$$

which means that $I - T$ is monotone. Similarly, the $\rho$-strong monotonicity of $F$ is equivalent to that of $I - T$.

For fixed points in the lower-level problem, consider the bifunction

$$F(u,v) = \langle u - T(u), v - u \rangle,$$

where $T : C \to C$ is a nonexpansive mapping. Assume that

$$\text{Fix}(T) := \{x \in C : T(x) = x\} \neq \emptyset.$$

It is clear that $F$ satisfies conditions $(C1)$–$(C3)$ in $(H1)$. We now verify $(C4)$. Since $T$ is nonexpansive, for all $u, v \in C$ we have

$$\langle T(u) - T(v), u - v \rangle \leq \|T(u) - T(v)\| \, \|u - v\| \leq \|u - v\|^2.$$

Thus,

$$F(u,v) + F(v,u) = \langle T(u) - T(v), u - v \rangle - \|u - v\|^2 \leq 0,$$

showing that $F$ is monotone.

When $C = \mathcal{H}$, we deduce that $x$ is an equilibrium point of $F$ if and only if $x \in \text{Fix}(T)$. When $C \neq \mathcal{H}$, $x \in \text{Fix}(T)$ implies $F(x, y) = 0$ for all $y \in C$. Hence, for every $p \in \mathcal{H}$,

$$\mathcal{F}_F(x, p) := \sup_{v \in C} \{\langle p, v \rangle + F(v, x)\} = \sigma_C(p),$$

and consequently, condition $(H5)$ reduces to

$$\limsup_{k \to +\infty} \lambda_k r_k \left[ \sigma_C \left( \tfrac{2p}{\lambda_k} \right) - \sigma_{\text{Fix}(T)} \left( \tfrac{2p}{\lambda_k} \right) \right] \leq 0, \quad \text{for all } x \in \text{Fix}(T), \ p \in \mathcal{N}_{\text{Fix}(T)}(x). \tag{30}$$

Note that

$$\sigma_{\text{Fix}(T)}\left( \tfrac{2p}{\beta_k} \right) = \tfrac{2}{\beta_k} \langle p, x \rangle, \quad p \in \mathcal{N}_{\text{Fix}(T)}(x),$$

and

$$\sigma_C\left( \tfrac{2p}{\beta_k} \right) = \sup_{u \in C} \langle p, u \rangle \leq \tfrac{2}{\beta_k} \sigma_C(p).$$

Thus, for all $x \in \text{Fix}(T)$ and $p \in \mathcal{N}_{\text{Fix}(T)}(x)$,

$$\lambda_k r_k \left[ \sigma_C \left( \tfrac{2p}{\lambda_k} \right) - \sigma_{\text{Fix}(T)} \left( \tfrac{2p}{\lambda_k} \right) \right] \leq 2 r_k \left[ \sigma_C(p) - \langle p, x \rangle \right].$$

If, moreover, $\sigma_C(p) < \infty$ for each $p \in \mathcal{N}_{\text{Fix}(T)}(x)$, then condition (30) is automatically satisfied whenever $\lim_{k \to +\infty} r_k = 0$.

The last requirement on the support functional $\sigma_C$ means that the barrier cone

$$\mathbf{b}(C) := \{p \in \mathcal{H} : \sigma_C(p) < \infty\}$$

contains $\mathcal{N}_{\text{Fix}(T)}(x)$ for every $x \in \text{Fix}(T)$. Using [41, Proposition 3.10], the polar cone of $\mathbf{b}(C)$ coincides with the asymptotic cone

$$S_C(x_0) := \bigcap_{t > 0} t(C - x_0), \quad x_0 \in C.$$

Similarly, by [41, Proposition 4.5], the tangent and normal cones

$$\mathcal{T}_{\text{Fix}(T)}(x) := \overline{\bigcup_{h>0} \tfrac{1}{h}\big(\text{Fix}(T) - x\big)}, \quad \mathcal{N}_{\text{Fix}(T)}(x),$$

are polar to each other, for every $x \in \text{Fix}(T)$. By polarity, the inclusion $\mathcal{N}_{\text{Fix}(T)}(x) \subset \mathcal{T}_{\text{Fix}(T)}(x)$, for every $x \in \text{Fix}(T)$, is equivalent to

$$S_C(x_0) \subset \mathcal{T}_{\text{Fix}(T)}(x), \quad \text{for all } x \in \text{Fix}(T), \text{ for all } x_0 \in C.$$

## 7 Conclusion and Perspectives

This work complements and extends both classical studies and recent investigations on prox-penalization, inertial proximal methods, and splitting algorithms for hierarchical equilibrium problems [17, 25–27]. In this chapter, we introduced the Relaxed Inertial Proximal Splitting Algorithm (RIPSA) for addressing such problems. Our convergence analysis began with weak ergodic and weak convergence of the generated iterates, established without assuming the presence of a Browder-Halpern contraction. We then proved strong convergence under the additional assumption of strong monotonicity on the upper-level bifunction. To mitigate the restrictiveness of this condition, we further incorporated a Browder-Halpern contraction factor into the algorithm, thereby recovering strong convergence results under weaker assumptions. Looking ahead, several research directions emerge. A natural extension is to investigate whether the same convergence techniques - especially those leading to strong convergence - can be adapted to Forward-Backward type algorithms, which are particularly effective for problems involving the minimization of a sum of two functions, one smooth and the other non-smooth (e.g., classical norm-based regularizers). Another promising avenue is to explore convergence under weaker assumptions, such as relaxed monotonicity conditions on bifunctions, which would enhance the applicability of the method in more realistic and complex settings.

## Acknowledgments

The authors are grateful to the referee for valuable comments, which have improved the presentation of the initially submitted version of this chapter.